\documentclass[12pt]{article}

\usepackage{amsmath, amsfonts, amsthm, amssymb, amscd}

\newtheorem{theorem}{Theorem}
\newtheorem{lemma}{Lemma}

\newcommand{\C}{\mathbb{C}}
\newcommand{\Q}{\mathbb{Q}}

\newcommand{\Hn}{{\mathcal H}_n}

\begin{document}

\title{Hodge structure on the Cohomology of \\
the Moduli space of Higgs bundles}
\author{Mridul Mehta \\
Department of Mathematics, University of Chicago}

\maketitle

\begin{abstract}
Let $C$ be a smooth projective curve over $\C$ of genus $g \ge 2$. 
Let $\Hn$ be the moduli space of stable Higgs bundles of rank $r$ and
degree $d$ on $C$ with values in $K \otimes \mathcal{O}(np)$, where
$K$ is the canonical bundle on $C$, $p$ any marked point in $C$, and 
$(r,d) = 1$. We prove that the natural Hodge structure on $H^k(\Hn, \Q)$
is pure of weight $k$.
\end{abstract}

An object that has been studied in great depth for many years now is the 
moduli space of flat unitary connections on a compact Riemann surface. Due
to the famous theorem of Narasimhan and Seshadri \cite{NS}, this space has 
been investigated by algebraic geometers under the guise of the moduli space
of vector bundles on a compact Riemann surface. More recently, several
people have focussed their attention on the moduli space of {\it all}
flat connections (as opposed to just unitary ones). There is a similar
correspondence theorem here as well, which identifies (only topologically
as in the former) this space with the moduli space of Higgs bundles on the
compact surface.

We know that under suitable conditions the moduli space of vector bundles
is a fine moduli space, so there exists a (holomorphic) universal bundle on this 
space. It was shown by Atiyah and Bott \cite{AB} that the Chern classes of this 
bundle can be decomposed appropriately into components known as universal 
classes, which then generate the cohomology ring of the moduli space. Recently, 
a similar result was proven in the case of Higgs bundles by Thaddeus and Hausel
\cite{TH1} in rank 2, and by Markman \cite{Mar} in general rank. In this note, 
we use their results to compute the Hodge structure on the cohomology of the 
moduli space of Higgs bundles on a smooth projective curve. The result is 
somewhat surprising since the moduli space is non-compact yet all of its 
cohomology groups have pure weight.

{\it Outline.}  This paper is divided into three parts. In section 1, we give a 
brief introduction to the moduli space of Higgs bundles and recall some of the 
results and constructions related to these objects. In section 2, we state the 
main result. Finally, section 3 contains the proof of the result.

{\it Notation and Conventions.} 
Throughout the paper, $C$ will denote the compact Riemann surface or smooth
projective curve, of genus $g \ge 2$. We shall use $\Hn$ to denote moduli 
spaces over $C$ of Higgs bundles $(E,\phi)$ of rank $r$ and degree $d$ with 
values in $K(n) = K \otimes {\mathcal O}(np)$, where $p$ is a marked point in 
$C$, and $n$ is any non-negative integer. We will always assume that $r$ and $d$ 
are coprime. All cohomology is with rational coefficients unless otherwise 
stated.

{\it Acknowledgements.}  I wish to thank my advisor Kevin Corlette, who first
suggested this problem to me. His guidance and support throughout have been
invaluable. I would also like to thank Vladimir Baranovsky for several
discussions which were immensely inspiring and helped me better understand
the mathematics related to this subject, and Madhav Nori for his valuable
insights and comments.

\section{Higgs bundles}

Higgs bundles were first studied in great detail by Hitchin and Simpson 
\cite{Hi, Si1}.
Let $C$ be a smooth complex projective curve of genus $g \ge 2$.
Let $L$ be any holomorphic line bundle on $C$. A {\it Higgs bundle} on $C$ with 
values in $L$ is a pair $(E, \phi)$, where $E$ is a holomorphic vector bundle 
on $C$, and $\phi$, often called a Higgs field, is an element of
$H^0(C,$ End$(E)\otimes L)$. We say that a Higgs bundle is {\it semi-stable} 
if for any $\phi$-invariant subbundle $F \subset E$, $\textrm{deg} F/
\textrm{rk} F \le \textrm{deg} E/\textrm{rk} E$. 
We say that the Higgs bundle is {\it stable} if this inequality
is strict. We will denote by $r$ the rank of $E$ and by $d$ the degree of $E$.

We now specialize to the case when $L = K(n) = K \otimes \mathcal{O}(np)$, 
where $K$ is the canonical bundle on $C$, $p$ is a marked point in $C$, and 
$n \ge 0$. Then the following is known.

\vspace{0.05 in}

\noindent
{\bf Theorem (Hitchin, Simpson, Nitsure \cite{Hi, Si1, Ni}).} \;
{\it For fixed rank $r$, and degree $d$ coprime to $r$, and for any
$n \ge 0$, there exists a moduli space $\Hn$ of stable Higgs bundles with
values in $K(n)$, which is a smooth quasi-projective variety of dimension
$r^2(2g - 2 + n) + 2$.}

\vspace{0.05 in}

When $n = 0$, the space $\mathcal{H}_0$ of Higgs bundles is related to the 
space of flat connections (or connections of constant central curvature) by
the following result.

\vspace{0.05 in}

\noindent
{\bf Theorem (Hitchin, Corlette, Donaldson \cite{Hi, Co, Do}).}
{\it Let $r$ and $d$ be coprime. Denote by $\mathcal{H}$ the space of
GL($r,\C$)-connections on $C$ of constant curvature $di \omega I$ up to
gauge equivalence, where 
$\omega$ is a 2-form on $C$ chosen so that $\int_C \omega = 2\pi/r$, and $I$
is the $r \times r$ identity matrix. Then $\mathcal{H}$ is a smooth variety
diffeomorphic to the moduli space $\mathcal{H}_0$.}

\vspace{0.05 in}

This equivalence is only topological, since the complex structures
on the two moduli spaces are different.

We now focus our attention on $\Hn$. This space can be constructed 
gauge-theoretically in the following manner (see \cite{TH1} 
for details). 
Let $\mathcal{V}$ be a Hermitian vector bundle on $C$ of rank $r$ and 
degree $d$. We denote by $\mathcal{A}$ the complex affine space of
holomorphic structures on $\mathcal{V}$. Let $\tilde{\Omega}^{p,q}$ be the
Sobolev completion of the space of differential forms of Sobolev class
$L^2_{k-p-q}$. For $n \ge 0$, we can then define 
a map
$$ \bar{\partial}: \mathcal{A} \times \tilde{\Omega}^{0,0}(\textrm{End} 
\;\mathcal{V}
\otimes K(n)) \rightarrow \tilde{\Omega}^{0,1}(\textrm{End}
\;\mathcal{V} \otimes
K(n))$$
by $\bar{\partial}(E, \phi) = \bar{\partial}_E \phi$, and let
$\mathcal{B}_n = \bar{\partial}^{-1}(0)$. Let $\mathcal{G}$ be the
gauge group of all complex automorphisms of $\mathcal{V}$, and 
$\overline{\mathcal{G}}$ be the quotient of $\mathcal{G}$ by the central 
subgroup
$\C^{\times} \subset \mathcal{G}$. Then we know that $\mathcal{G}$ acts 
on $\mathcal{A}$ and
$\tilde{\Omega}^{1,0}(\textrm{End} \; \mathcal{V} \otimes K(n))$ 
thus inducing
an action on $\mathcal{B}_n$. Finally, if we denote by $\mathcal{B}^s_n
\subset \mathcal{B}_n$ the open subset of stable Higgs bundles, then we have
$\Hn = \mathcal{B}^s_n/\overline{\mathcal{G}}$. 

We also know that for all $n \ge 0$, $\Hn \times C$ carries on it a
universal family $(\mathbb{E}_n,\Phi_n)$. Although the universal bundles 
$\mathbb{E}_n$ are not canonical, the projective bundles 
$\mathbb{P}(\mathbb{E}_n)$ are, and we recall briefly one method of 
constructing these bundles.

We start with the tautological rank $r$ vector bundle 
$\mathbb{E}_\mathcal{A}$ on $\mathcal{A} \times C$ with the constant
scalars $\C^{\times} \subset \mathcal{G}$ acting trivially on the base
and as scalars in the fibre of $\mathbb{E}_\mathcal{A}$. Using the 
projection map $\pi:
\mathcal{B}^s_n \rightarrow \mathcal{A}$, we pull back the bundle
$\mathbb{E}_\mathcal{A}$ to get a bundle we will
call $\mathbb{E}_{\mathcal{B},n}$ on $\mathcal{B}^s_n \times C$. Next, 
we notice that the group $\mathcal{G}$ acts on this bundle, thus inducing
an action on the projective bundle 
$\mathbb{P}(\mathbb{E}_{\mathcal{B},n})$ with the
constant scalars $\C^{\times} \subset \mathcal{G}$ acting trivially on the
fibres. Hence the bundle $\mathbb{P}(\mathbb{E}_{\mathcal{B},n})$ 
is a $\overline{\mathcal{G}}$-equivariant
bundle. Since $\overline{\mathcal{G}}$ also acts freely on $\mathcal{B}^s_n$,
we see that $\mathbb{P}(\mathbb{E}_{\mathcal{B},n})$ descends to a 
$\mathbb{P}^{r-1}$-bundle
$\mathbb{P}(\mathbb{E}_n)$ on $\Hn \times C$. Any lift of
this bundle is a universal bundle on $\Hn \times C$. In general of course
such lifts don't always exist. However, for the spaces $\Hn$, one can
always find a lift and thereby construct a universal bundle.
The crucial point to note here is that the projective bundle 
$\mathbb{P}(\mathbb{E}_n)$ on $\Hn \times C$ is constructed by pulling back
the bundle $\mathbb{E}_\mathcal{A}$ which is
independent of the parameter $n$. This allows us to do the following.

Thaddeus and Hausel \cite{TH1} showed that the spaces $\Hn$ constructed in the 
above manner form a resolution tower, that is there are natural
inclusions $\Hn \hookrightarrow \mathcal{H}_{n+1}$, which allow us to 
construct the direct limit space which we shall denote by 
$\mathcal{H}_{\infty}$. They also showed that the $\mathcal{B}_n$'s also have
natural inclusions which allows us to construct the direct limit 
$\mathcal{B}_{\infty}$ and that $\mathcal{B}^s_{\infty}$ is a principal 
$\overline{\mathcal{G}}$-bundle on $\mathcal{H}_\infty$ with total space
contractible. Since the projections $\mathcal{B}_n \rightarrow
\mathcal{A}$ commute with the inclusions $\mathcal{B}_n \hookrightarrow
\mathcal{B}_{n+1}$, and the construction of the bundles 
$\mathbb{P}(\mathbb{E}_n)$ did not depend on $n$, we can now take the
direct limit of these bundles and thus obtain a projective bundle 
on $\mathcal{H}_\infty$ which
we shall denote by $\mathbb{P}(\mathbb{E}_{lim})$. We also note that due
to the way it is constructed, $\mathbb{P}(\mathbb{E}_{lim})$ is the
canonical bundle $\mathbb{P}(\textrm{\bf E})$ where {\bf E} is any universal
bundle on $\mathcal{H}_\infty$. 

Finally, we will denote by ${c}^n_1, \ldots, {c}^n_r$
the Chern classes of the universal bundle ${\mathbb E}_n$ on $\Hn$,
and by $\overline{c}^n_2, \ldots, \overline{c}^n_r$ the Chern classes of the
projective bundle $\mathbb{P}({\mathbb E}_n)$ on $\Hn$. The $\overline{c}^n_k$
are elements of rational cohomology, and may be thought of as Chern
classes of the universal
bundle $\mathbb{E}$ twisted by a formal $r$-th root of 
$\Lambda^r\mathbb{E}^*$, so that the first Chern class vanishes. 
Similarly, $\overline{c}_2, \ldots, \overline{c}_r$ will denote the Chern 
classes of 
the projective bundle $\mathbb{P}({\mathbb E}_{lim})$ on 
${\mathcal H}_{\infty}$. 
One consequence of the above construction is that the Chern classes of 
$\mathbb{P}(\mathbb{E}_{lim})$ restrict to their counterparts on $\Hn$ for
$n \ge 0$. 
The bundle $\mathbb{P}(\mathbb{E}_{lim})$ will play a key role in our
proof in section 3.

\section{Statement of the Result}

\begin{theorem}
The natural Hodge structure on $H^k(\Hn,\Q)$ is pure of weight $k$.
\end{theorem}

This is obviously known to be true in the case of smooth projective varieties
over $\C$. However, if the variety is singular or non-compact, we typically
expect the natural Hodge structure on $H^k$ to have some {\it mixing}, i.e., 
we expect the $k$-th cohomology group to contain elements of weight other 
than $k$. 
So it is interesting and perhaps a little suprising that although $\Hn$ is 
non-compact, there is indeed no such {\it mixing} in its cohomology.

A simple consequence of the above result is the following. If $X$ is any
smooth compactification of $\Hn$, then the natural maps
$H^k(X, \Q) \rightarrow H^k(\Hn, \Q)$ are surjective. This follows by
first considering the case when the compactification divisor is a normal 
crossings divisor, and then using Hironaka's resolution method to always
reduce to this case.

The proof of the above theorem occupies the rest of the paper. We use
the approach of Thaddeus and Hausel \cite{TH1} to get a grip on the spaces 
$\Hn$. More specifically,
we prove first that the $k$-th cohomology space $H^k(\mathcal{H}_\infty, \Q)$
of the direct limit $\mathcal{H}_\infty$ of the spaces $\Hn$ carries a 
natural Hodge structure. We then show that this in fact is pure, using the
fact that the cohomology is generated by the universal classes of the
direct limit bundle $\mathbb{P}(\mathbb{E}_{lim})$. This finally allows us 
to use the
surjective map from $H^*(\mathcal{H}_\infty) \rightarrow H^*(\Hn)$ to complete
the proof. 

\section{Proof}

\begin{lemma}
The cohomology ring $H^*({\mathcal H}_{\infty}, \Q)$ carries a natural mixed 
Hodge structure (MHS).
\end{lemma}

\begin{proof}
Since each $\Hn$ is a smooth quasi-projective variety, it carries 
a natural MHS \cite{De}.
Let $W_{n,*}$ and $F^{n,*}$ denote the weight and Hodge filtrations 
respectively on $H^*(\Hn,\Q)$.
Define filtrations $W_*$ and $F^*$ on $H^*({\mathcal H}_{\infty}, \Q)$ by 
setting $W_* = \underset{n}{\varprojlim} W_{n,*}$ and 
$F^* = \underset{n}{\varprojlim} F^{n,*}$. It is not difficult to 
verify that these filtrations satisfy the following properties:
\begin{enumerate}
\item[(i)] $W_*$ defines a weight filtration on $H^*({\mathcal H}_{\infty},
\Q)$
\item[(ii)] $F^*$ defines a Hodge filtration on $H^*({\mathcal H}_{\infty},
\Q) \otimes \C$
\item[(iii)] $F^*$ induces a $\Q$-Hodge filtration of weight $k$ on each
$\textrm{Gr}_k^W(H^*({\mathcal H}_{\infty}, \Q)) = W_k/W_{k-1}$
\end{enumerate}
Consequently, we see that the filtrations $F^*$ and $W_*$ define a 
MHS on $H^*({\mathcal H}_{\infty}, \Q)$ in the sense of Deligne, although
this is an infinite-dimensional vector space.
\end{proof}

\begin{lemma}
Let $\Gamma_{k,i}$ denote the composition of the following maps.

$$\begin{CD}
  H^i(C \times \mathcal H_\infty, \Q)(-k+1)	@>\cup\overline{c}_k>>	
H^{i+2k}(C \times \mathcal H_\infty, \Q)(1) 	\\
  @AAp_1^*A							@VV\int_CV		
			\\
  H^i(C, \Q)(-k+1)				&&			
H^{i+2k-2}(\mathcal H_\infty, \Q)
\end{CD}$$

\vspace{0.13in}

\noindent
Then the cohomology $H^*({\mathcal H}_{\infty}, \Q)$ is generated as a ring
by the images of the maps $\Gamma_{k,i}$, $2 \le k \le r$, $0 \le i \le 2$.
\end{lemma}

\begin{proof}
See Thaddeus and Hausel \cite{TH1}, 10.1.
\end{proof}

Now we focus our attention on $\Hn$. Since $\Hn$ is smooth quasi-projective, 
we may assume that
$\Hn = X_n \backslash D_n$ where $X_n$ is smooth projective and $D_n$ is a 
normal crossings divisor \cite{Hi}. Let $\Omega^p(*D_n)$ denote the sheaf 
on $X_n$ of meromorphic $p$-forms that are holomorphic on $\Hn$ and have
poles of arbitrary (finite) order on $D_n$. Similarly, let 
$\mathcal{A}^p(*D_n)$ denote the sheaf on $X_n$ associated to the presheaf 
$U \rightarrow A^p(U \backslash U \cap D_n)$, where $A^p$ is the space of
smooth $p$-forms. Both of these fit into
complexes of sheaves $(\Omega^*(*D), d)$ and $(\mathcal{A}^*(*D), d)$
on $X_n$. Next, we let $\Omega^p(\textrm{log} D)$ be the subsheaf of
$\Omega^p(*D)$ generated locally by the holomorphic forms and the
logarithmic differentials $dz_i/z_i$, $1 \le i \le k$, where $D_n$ can
be written locally as $\{z_1 \cdots z_k = 0\}$. Intrinsically, if $f$
is a local defining equation of $D_n$, then $\Omega^p(\textrm{log} D)$
is given by those meromorphic $p$-forms $\psi$ such that both $f \psi$
and $f d\psi$ are holomorphic. We recall that both the inclusions
$\Omega^*(\textrm{log} D) \subset \mathcal{A}^p(*D_n)$ and
$\Omega^*(*D) \subset \mathcal{A}^p(*D_n)$ are quasi-isomorphisms. Since
the sheaves $\mathcal{A}^p(*D_n)$ are {\it fine}, 
$H^q(X_n, \mathcal{A}^p(*D_n)) = 0$ for $q > 0$, using the spectral 
sequence for hypercohomology, we see that 
$\mathbb{H}^*(X_n, \Omega^*(*D)) \cong
\mathbb{H}^*(X_n, \mathcal{A}^p(*D_n)) \cong 
H^*_d(H^0(X_n, \mathcal{A}^p(*D_n))) = H^*_{\textrm{DR}}(\Hn)$. 

Hence
we see that the cohomology of $\Hn$ can be computed by using meromorphic
forms on $X_n$ that are holomorphic on $\Hn$ and have poles along $D_n$.
This enables us to define the {\it weight} of a $(p,q)$ form $\gamma$
in the cohomology of $\Hn$ (see \cite{De} for more details). We say 
that the form $\gamma$ has weight 
$p+q+d$ where $d$ is the order of the pole along $D_n$ of the associated 
meromorphic form on $X_n$. In particular if $\gamma$ is the restriction
of a form holomorphic on all of $X_n$, then it has weight $p+q$ (i.e.,
$d = 0$).

Before proving the next lemma, we discuss briefly the notion of Chern classes
for coherent sheaves on a smooth projective variety. Let $X$ be a smooth
projective variety. Given a coherent sheaf 
$\mathcal{F}$ on $X$, we can consider $\mathcal{F}$ as an element of the 
Grothendieck group $G(X)$. We consider $K(X)$, defined in a similar manner
as $G(X)$ except using locally free sheaves, as the quotient
of the free abelian group generated by all locally free (coherent) sheaves
on $X$, by the subgroup generated by all expressions of the form
$\mathcal{E} - \mathcal{E}' - \mathcal{E}''$, whenever $0 \rightarrow
\mathcal{E'} \rightarrow \mathcal{E} \rightarrow \mathcal{E}'' \rightarrow 0$
is a short exactly sequence of locally free sheaves. 
We recall that the definition of Chern classes extends naturally to $K(X)$
(see for instance Hartshorne, {\it Algebraic Geometry}).
It can be shown using resolutions of coherent sheaves that
the natural map from $K(X) \rightarrow G(X)$ is an isomorphism. This allows
us to consider the image of the coherent sheaf $\mathcal{F}$ in $K(X)$,
and consequently define Chern classes of $\mathcal{F}$. This turns out to be
a good definition in that it satisfies all the defining properties of 
Chern classes.
We now prove the following.

\begin{lemma}
Let $\Gamma^n_{k,i}$ denote the composition of the following maps.

$$\begin{CD}
  H^i(C \times \Hn, \Q)(-k+1)	@>\cup \overline{c}^n_k>>	H^{i+2k}(C 
\times \Hn, \Q)(1) 	\\
  @AAp_1^*A					@VV\int_CV			
\\
  H^i(C, \Q)(-k+1)		&&			H^{i+2k-2}(\Hn, \Q)
\end{CD}$$

\vspace{0.13in}

\noindent
Then for each $k$, $i$, $\Gamma^n_{k,i}$ is a morphism of Hodge structures.
\end{lemma}

\begin{proof}
We note that $\Gamma^n_{k,i}$ is a composition of three maps: 
a pullback $p_1^*$,
cupping with the Chern class $\overline{c}^n_k$, followed by integration along 
the 
curve $C$. Pulling back a form
$\omega \in H^i(C,\Q)$ to $H^i(C \times \Hn, \Q)$  does not change the Hodge 
type or weight of
the form (since $p_1^*\omega$ will not have any poles along $D_n$), so this 
is  a morphism 
of Hodge structures. Next we observe that the universal bundle $\mathbb E_n$ 
on $\Hn$ can be
extended to a coherent sheaf $\mathcal E_n$ on $X_n$. Then 
the $k$-th Chern class of $\mathcal E_n$
restricts to give the Chern class $c^n_k$ of $\mathbb E_n$. Hence 
$\overline{c}^n_k$ (which we already
know to be of pure Hodge type $(k,k)$) does not have any poles along $D_n$ 
and therefore has weight
$2k$. Consequently, cupping the pullback $p_1^* \omega$ of any $\omega \in 
H^i(C,\Q)$ with $\overline{c}^n_k$
will affect the Hodge type as expected by $(k,k)$, but will also shift weight 
exactly by $2k$. As a result,
this is also a morphism of Hodge structures. Finally, integrating a form 
$p_1^*\omega \cup c^n_k$
along $C$ is also a morphism of Hodge structures since the 
Hodge type will shift down
by $(1,1)$ and weight will always shift exactly by 2 (since $p_1^*\omega 
\cup \overline{c}^n_k$ has no poles
along $D_n$).
\end{proof}

Thus each $\Gamma^n_{k,i}$ is a morphism of Hodge structures. We can now finish 
the proof.

\bigskip

\noindent
{\it Proof of Theorem 1.}
We now take the
limit of these maps as $n\rightarrow \infty$ to get a map from 
$H^i(C,\Q) \rightarrow H^{i+2k-2}(\mathcal H_\infty, \Q)$, which 
is exactly the map $\Gamma_{k,i}$ constructed before. Since the MHS on
$H^*(\mathcal H_\infty, \Q)$ was obtained as a limit of the MHS on 
$H^*(\Hn, \Q)$, $\Gamma_{k,i}$ is also a morphism of Hodge structures. 
By Lemma 2, we know that the image of $\Gamma_{k,i}$ 
generates the cohomology ring
$H^*(\mathcal H_\infty, \Q)$. Since $\Gamma_{k,i}$ is a morphism of Hodge 
structures, any $\omega \in 
H^k(\mathcal H_\infty, \Q)$ has weight $k$. Next we recall that the natural 
map $H^*(\mathcal H_\infty, \Q) \rightarrow H^*(\Hn, \Q)$ is 
surjective \cite{TH1} and because of the manner in which the MHS on 
$H^*(\mathcal H_\infty, \Q)$ was induced,
this is also a morphism of Hodge structures. Hence, we see that for any $k$, 
all elements of $H^k(\Hn, \Q)$ are of weight $k$. This completes the
proof. $\square$


\begin{thebibliography}{1}

\bibitem{AB} M.F. Atiyah and R. Bott, {\it The Yang-Mills equations over 
Riemann surfaces}, Philos. Trans. Roy. Soc. London Ser. A {\bf 308} (1982)
523-615.

\bibitem{Co} K. Corlette, {\it Flat $G$-bundles with canonical metrics},
J. Differential Geom. {\bf 28} (1988) 361-382.

\bibitem{De} P. Deligne, {\it Théorie de Hodge II}, Inst. Hautes 
Études Sci. Publ. Math. {\bf 40} (1971) 5-57.

\bibitem{Do} S.K. Donaldson, {\it Twisted harmonic maps and the self-duality
equations}, Proc. London. Math. Soc. {\bf 55} (1987) 127-131.

\bibitem{Ha1} T. Hausel, {\it Geometry of the moduli space of Higgs bundles},
Ph.D. thesis, University of Cambridge, 1998; available from the Front for
the Mathematics ArXiv AG/0107040.

\bibitem{Ha2} T. Hausel, {\it Compactification of moduli of Higgs bundles},
J. Reine Angew. Math. {\bf 503} (1998) 169-192.

\bibitem{Hi} H. Hironaka, {\it On resolution of singularities}, Proc. Int.
Congress Math., Stockholm (1962) 507-525.

\bibitem{Hit} N. Hitchin, {\it The self-duality equations on a Riemann
surface}, Proc. Lond. Math. Soc. {\bf 55} (1987) 59-126.

\bibitem{Mar} E. Markman, {\it Generators of the cohomology ring of moduli 
spaces of sheaves on symplectic surfaces}, J. Reine Angew. Math. {\bf 544} 
(2002) 61-82.

\bibitem{NS} M.S. Narasimhan and C.S. Seshadri, {\it Stable and unitary 
vector bundles on a compact Riemann surface}, Ann. of Math. (2)
{\bf 82} (1965) 540-567.

\bibitem{Ni} N. Nitsure, {\it Moduli Spaces of semistable pairs on a curve},
Proc. London Math. Soc. (3) {\bf 62} (1991) 275-300.

\bibitem{Si1} C.T. Simpson, {\it Higgs bundles and local systems}, Inst.
Hautes Etudes Sci. Publ. Math. {\bf 75} (1992) 5-95.

\bibitem{Si2} C.T. Simpson, {\it The ubiquity of variations of Hodge
structure}, Complex geometry and Lie theory (Sundance, UT, 1989), 329-348,
Proc. Symp. Pure Math. {\bf 53}, Amer. Math. Soc., 1991.

\bibitem{Si3} C.T. Simpson, {\it Moduli of representations of the 
fundamental group of a smooth projective variety I}, Inst. Hautes Sci. Publ.
Math. {\bf 79} (1994) 47-129.

\bibitem{Si4} C.T. Simpson, {\it Moduli of representations of the 
fundamental group of a smooth projective variety II}, Inst. Hautes Sci. Publ.
Math. {\bf 80} (1995) 5-79.

\bibitem{TH1} M. Thaddeus and T. Hausel, {\it Generators for the cohomology ring 
of the moduli space of rank 2 Higgs bundles}, available from the Front for
the Mathematics ArXiv AG/0003093.

\bibitem{TH2} M. Thaddeus and T. Hausel, {\it Relations in the cohomology ring 
of the moduli space of rank 2 Higgs bundles}, available from the Front for
the Mathematics ArXiv AG/0003094.

\end{thebibliography}
\end{document}